\documentclass[arxiv,numberwithout=true]{agn_article}
\PassOptionsToPackage{ngerman}{babel}
\usepackage[tikzexternalize=false]{agn_all}
\usepackage{oldgerm}
\usepackage{amsmath,mathtools,amssymb}
\usepackage{agn_macros}
\usepackage{float}
\usepackage{yfonts}
\usepackage{longtable}
\usepackage{enumitem}
\usepackage{perpage}
\usepackage[hyperindex]{hyperref}

\pgfplotsset{
every axis/.append style={width=0.6\textwidth, height=6cm},
every axis plot/.append style={no markers, thick},
label style={font=\small},
tick label style={font=\small},
legend pos={outer north east},
cycle list={red, green, cyan, yellow}
}
\DeclareMathOperator{\scalee}{e}
\newcommand{\force}{\widetilde{f}}
\newcommand{\const}{H}
\newcommand{\opoint}{p}%
\newcommand{\ipoint}{q}%
\newcommand{\deform}{\overline{F}}%
\newcommand{\deformT}{\smash{\overline{F}}\vphantom{F}^T}
\newcommand{\jacobi}{M}%
\newcommand{\work}{W}%
\newcommand{\hpressure}{\overline{\sigma}}
\renewcommand{\.}{\hspace*{0.07em}}
\newtheorem{postulat}{V\!}
\newtheorem{postulat2}{D\!}
\newtheorem{postulat3}{P\!}

\newcommand{\quoteref}[1]{\enquote{\emph{#1}}}
\newcommand{\quoteesc}[1]{{\rm[}#1{\rm]}}
\let\defaultfootnote\footnote%
\newcounter{transfootnote}
\MakePerPage{transfootnote}
\newcounter{footnotetemp}
\newcommand{\fnsymbolfootnote}[1]{%
	\setcounter{footnotetemp}{\value{footnote}}%
	\let\thefootnotedefault\thefootnote%
	\renewcommand{\thefootnote}{\realfnsymbol{transfootnote}}%
		\ifnum\value{transfootnote}=7%
			\setcounter{transfootnote}{0}%
		\fi%
		\refstepcounter{transfootnote}%
		\defaultfootnote{#1}%
	\let\thefootnote\thefootnotedefault%
	\setcounter{footnote}{\value{footnotetemp}}%
}
\newcommand{\transcomment}[1]{%
	\fnsymbolfootnote{Translators' remark: #1}%
}

\begin{document}
\begin{refsegment}
\selectlanguage{\english}

\title{The axiomatic introduction of arbitrary strain tensors by Hans Richter - a commented translation of\\ \enquote{Strain tensor, strain deviator and stress tensor for finite deformations}}
\author{%
	Patrizio Neff\thanks{%
		Patrizio Neff,\quad Head of Lehrstuhl f\"{u}r Nichtlineare Analysis und Modellierung, Fakult\"{a}t f\"{u}r	Mathematik, Universit\"{a}t Duisburg-Essen, Thea-Leymann Str. 9, 45127 Essen, Germany, email: patrizio.neff@uni-due.de%
	}\quad\ and \quad%
	Kai Graban\thanks{%
		Kai Graban, \quad Fakult\"{a}t f\"{u}r Mathematik, Universit\"{a}t Duisburg-Essen, Thea-Leymann Str. 9, 45127 Essen%
	}\quad\ and\quad%
	Eva Schweickert\thanks{%
		Eva Schweickert,\quad Lehrstuhl f\"{u}r Nichtlineare Analysis und Modellierung, Fakult\"{a}t f\"{u}r Mathematik, Universit\"{a}t Duisburg-Essen, Thea-Leymann Str. 9, 45127 Essen, Germany; email: eva.schweickert@stud.uni-due.de%
	}\quad\ and\quad%
	Robert J.\ Martin\thanks{%
		Robert J.\ Martin,\quad Lehrstuhl f\"{u}r Nichtlineare Analysis und Modellierung, Fakult\"{a}t f\"{u}r Mathematik, Universit\"{a}t Duisburg-Essen, Thea-Leymann Str. 9, 45127 Essen, Germany; email: robert.martin@uni-due.de%
	}%
}
\date{\today}
\maketitle
\begin{abstract}
	We provide a faithful translation of Hans Richter's important 1949 paper \enquote{Verzerrungstensor, Verzerrungsdeviator und Spannungstensor bei endlichen Formänderungen} from its original German version into English, complemented by an introduction summarizing Richter's achievements.
\end{abstract}
\\[.65em]
\let\footnote\fnsymbolfootnote%
\section*{Introduction}
In this paper, we continue our efforts to translate Hans Richter's early work on nonlinear elasticity theory (cf.\ \cite{agn_graban2019richter}). Richter's second article in the field, entitled \enquote{Verzerrungstensor, Verzerrungsdeviator und Spannungstensor bei endlichen Formänderungen} (\enquote{Strain tensor, strain deviator and stress tensor for finite deformations}) \cite{richter1949verzerrung}, was published in \emph{Zeitschrift für Angewandte Mathematik und Mechanik} in 1949 and concerns the axiomatic foundations of nonlinear elasticity. More precisely, Richter is concerned with introducing deductively a family of strain tensors for which he lays down an axiomatic structure.

In order to provide the context for Richter's work, we briefly recapitulate what can be said, and what is generally accepted, about strain tensors, following Truesdell and his school after 1955. The concept of \emph{strain} is of fundamental importance in continuum mechanics. In linearized elasticity, it is assumed that the Cauchy stress tensor $\sigma$ is a linear function of the symmetric infinitesimal strain tensor
\[
	\eps = \sym \grad u = \sym(\grad \varphi - \id) = \sym(F-\id)\,,
\]
where $\varphi\col\Omega\to\R^n$ is the deformation of an elastic body with a given reference configuration $\Omega\subset\R^n$, $\varphi(x) = x+u(x)$ with the displacement $u$, $F=\grad\varphi$ is the deformation gradient, $\sym\grad u = \frac12(\grad u + (\grad u)^T)$ is the symmetric part of the displacement gradient $\grad u$ and $\id$ is the identity tensor.
In geometrically nonlinear elasticity, on the other hand, a vast number of different \enquote{strains} have been employed in the past in order to conveniently express nonlinear constitutive relations. In particular, it is common practice to choose a stress-strain pair such that a given constitutive law can be expressed in terms of a linear relation between stress and strain \cite{batra1998linear,batra2001comparison,bertram2007rank}.%
\footnote{%
	Cf.\ Truesdell and Noll \cite[p.~347]{truesdell65}: \quoteref{Various authors \quoteesc{\ldots} have suggested that we should select the strain \quoteesc{tensor} afresh for each material in order to get a simple form of constitutive equation. \quoteesc{\ldots} \emph{Every} invertible stress relation $T=f(B)$ for an isotropic elastic material is linear, trivially, in an appropriately defined, particular strain \quoteesc{tensor $f(B)$}.}%
}
In these cases, the strain tensor is generally a nonlinear function of the deformation gradient.

Although the specific definition of what exactly the term \enquote{strain} encompasses varies throughout the literature, it is commonly assumed \cite[p.\ 230]{Hill68} (cf.\ \cite{hill1970, hill1978, bertram2008elasticity, norris2008higherDerivatives}) that a (\emph{spatial} or \emph{Eulerian}) strain takes the form of a \emph{primary matrix function} of the left Biot-stretch tensor $V=\sqrt{FF^T}$ of the deformation gradient $F\in\GLpn$, i.e.\ an isotropic tensor function $E\col\PSymn\to\Symn$ from the set of positive definite tensors to the set of symmetric tensors of the form%
\footnote{%
	Note that more general definitions can be found in the literature as well \cite{agn_neff2015geometry,truesdell60}; for example, Truesdell and Toupin \cite[p.~268]{truesdell60} consider \quoteref{any uniquely invertible isotropic second order tensor function of \quoteesc{the left Cauchy-Green deformation tensor $B=FF^T$}} to be a strain tensor.
}
\begin{equation}
\label{eq:primaryMatrixFunctionDefinition}
	E(V) = \sum_{i=1}^n \scalee(\lambda_i) \cdot e_i\otimes e_i \quad\text{for}\quad V = \sum_{i=1}^n \lambda_i \cdot e_i\otimes e_i
\end{equation}
with a strictly monotone \emph{scale function} $\scalee\col(0,\infty)\to\R$, where $\otimes$ denotes the tensor product, $\lambda_i$ are the eigenvalues and $e_i$ are the eigenvectors of $V$. In addition, the normalization requirements $\scalee(1)=0$ and $\scalee'(1)=1$ are typically required to hold as well, with the former ensuring that the strain vanishes if and only if the deformation gradient describes a pure rotation, i.e.\ if and only if $F\in\SO(n)$, where $\SO(n)=\{Q\in\GLn \setvert Q^TQ=\id,\, \det Q = 1\}$ denotes the special orthogonal group. This property, in turn, ensures that the only strain-free deformations are rigid body movements \cite{Neff_curl06}.

\subsection*{Richter's general definition of strain}

We now turn to Richter's original development, which precedes the work of Truesdell.
Based on the polar decomposition $F=V\.R=R\.U$ with $R\in\SO(3)$ and $U,V\in\PSym(3)$ of the deformation gradient $F\in\GLp(3)$ as well as a certain notion of superposition (which is described in more detail in the following section), Richter arrives at a fully general definition of
Eulerian as well as Lagrangian strain tensors. Expressed in terms of the \emph{principal matrix logarithm} $\log\col\PSymn\to\Symn$ on the set $\PSymn$ of positive definite symmetric matrices, Richter's definition is given by
\begin{subequations}
\label{straintensordef}
	\begin{alignat}{2}
		E(F) &= \ftilde(\log V)\in\Sym(3)\qquad\qquad &&\text{(Eulerian strains),}\label{eq:straintensordefEulerian}\\
		\widehat{E}(F) &= \ftilde(\log U)\in\Sym(3)\qquad\qquad &&\text{(Lagrangian strains),}\label{eq:straintensordefLagrangian}
	\end{alignat}
\end{subequations}
where $\ftilde\colon\Sym(3)\to\Sym(3)$ is any differentiable and invertible (i.e.\ injective) primary matrix function\footnote{%
	Here and throughout, we will identify the primary matrix function with its associated scale function and write, for example, $\ftilde(V)=\sum_{i=1}^n \ftilde(\lambda_i) \cdot e_i\otimes e_i$.%
}
of the form \eqref{eq:primaryMatrixFunctionDefinition} with $\ftilde(0)=0$ and $\ftilde'(0)=1$.
In particular, due to the invertibility of the principal matrix logarithm, Richter's definition is indeed equivalent to the contemporary definition \eqref{eq:primaryMatrixFunctionDefinition} of a general strain tensor; note that since $\scalee(1)=\ftilde(0)$ and $\scalee'(1)=\ftilde'(0)$ for $\ftilde=\scalee\circ\exp$ and $\scalee=\ftilde\circ\log$, the stated normalization requirements are equivalent as well.

Similar to Richter, we will mostly focus on the Eulerian family \eqref{eq:straintensordefEulerian} in the following; analogous considerations can of course be applied to the Lagrangian family as well.
First, note that the invertibility of $\ftilde$ implies the equivalence
\begin{align*}
	E(F)=0
	\quad\iff\quad \log V=0
	\quad\iff\quad V=\id
	\quad\iff\quad F\in\SO(3)\,,
\end{align*}
thus $E(\grad\varphi)\equiv0$ if and only if $\varphi$ is a rigid body movement \cite{agn_munch2008curl}.
Furthermore, Richter's definitions \eqref{straintensordef} naturally contain a number of commonly employed strains, including the material and spatial \emph{Hencky strain tensors} \cite{Hencky1928,Hencky1929,hencky1929super,hencky1931,agn_neff2015geometry,agn_neff2015exponentiatedI,agn_neff2015exponentiatedII,agn_neff2014exponentiatedIII,agn_neff2014rediscovering}
\begin{equation}
	E_0 =\log V = \log(\sqrt{FF^T})\,,
	\qquad
	\widehat{E}_0 =\log U = \log(\sqrt{F^TF})\,,
\end{equation}
which are often been considered to be the \emph{natural} or \emph{true} strains in nonlinear elasticity \cite{tarantola2009stress,Tarantola06,freed1995natural,hanin1956isotropic},
as well as the Seth-Hill \cite{seth1935,Hill68,seth1961generalized} and Doyle-Ericksen \cite{doyle1956nonlinear} strain tensor families
\begin{align}
	E_{(m)}=\frac{1}{2m}(V^{2m}-\id)=\frac{1}{2m}(B^m-\id)\,,
	\qquad\qquad
	\widehat{E}_{(m)}=\frac{1}{2m}(U^{2m}-\id)=\frac{1}{2m}(C^m-\id)\,.
\end{align}
However, Richter's definition \eqref{straintensordef} is significantly more general and includes, for example, the \emph{Ba\v zant strain tensor} \cite{bazant1995}, given by $\frac12(V-V\inv)$; note that for $\ftilde(\lambda)=\frac12(e^\lambda-e^{-\lambda})$ or, equivalently, $\ftilde\inv(x)=\log(x+\sqrt{x^2+1})$,
\begin{align}
	\ftilde(\log V) = \frac12\.(\exp(\log V) - \exp(\log V)\inv) = \frac12\.(V-V\inv)\,.
\end{align}
Another example is the (Eulerian) Almansi strain tensor \cite{almansi1911sulle}, attributed to Trefftz in a review of Richter's article by Moufang, which is given by $T=\frac{1}{2}(\id-B\inv)$ with $B=V^2$ and corresponds to the choice $\ftilde(\lambda)=\frac12(1-e^{-2\lambda})$ for the transition function $\ftilde$ in \eqref{eq:straintensordefEulerian}.
\begin{figure}[H]
\tikzremake
	\centering
	\begin{minipage}[c]{0.45\linewidth}
	\centering
	\begin{tikzpicture}
	\begin{axis}[axis x line=center,
		axis y line=center,
		xtick={-4,-3,-2,-1,0,1,2,3,4},
		ytick={-4,-3,-2,-1,0,1,2,3,4},
		xticklabels={},
		yticklabels={},
		xlabel={$\lambda$},
		ylabel={$\ftilde(\lambda)$},
		xlabel style={below right},
		ylabel style={above left},
		xmin=-5,
		xmax=5,
		ymin=-5,
		ymax=5,
		samples=420,
		width=.931\linewidth,
		height=.6\linewidth,
		domain=-2:4.2
	]
		\addplot{0.5*(1-exp(-2*x))};
	\end{axis}
	\end{tikzpicture}
	\end{minipage}	
	\caption{Transition function $\ftilde$ for the Almansi strain tensor $T=\frac{1}{2}(\id-B\inv)$.}
	\label{fig:transition}
\end{figure}
Observe that Richter's strain tensors are isomorphic to each other%
\footnote{%
	Cf.\ Truesdell and Toupin \cite[p.~268]{truesdell60}: \quoteref{\ldots any \quoteesc{tensor} sufficient to determine the directions of the principal axes of strain and the magnitude of the principal stretches may be employed and is fully general}. Truesdell and Noll \cite[p.~348]{truesdell65} also argue that there \quoteref{is no basis in experiment or logic for supposing nature prefers one strain \quoteesc{tensor} to another}.
}
in the sense that for any pair $E_1$, $E_2$ of strain tensors in the family \eqref{eq:straintensordefEulerian}, there exists an invertible, isotropic mapping $\zeta\colon\Sym(3)\to\Sym(3)$ such that
\begin{align}
	E_1=\zeta(E_2)\,;
\end{align}
since $E_1=\ftilde_1(\log V)$ and $E_2=\ftilde_2(\log V)$ for suitable invertible functions $\ftilde_1,\ftilde_2$, it suffices to choose $\zeta=\ftilde_1\circ\ftilde_2\inv$.

We also note that a strain tensor $E$ of the form \eqref{eq:straintensordefEulerian} is \emph{tension-compression symmetric}, i.e.\ satisfies $E(V\inv)=-E(V)$, if and only if $\ftilde$ is odd, i.e.\ if $\ftilde(\lambda)=-\ftilde(-\lambda)$.

\subsection*{Richter's superposition principle}
Richter obtains his general definition \eqref{straintensordef} deductively from three axioms. Most importantly, he assumes that any strain tensor satisfies a \emph{superposition principle} (postulate \textbf{V3}) in the case of \emph{coaxial} stretches. More specifically, for $V_1,V_2\in\PSym(3)$ such that $V_1V_2=V_2V_1$, let $E_1=E(V_1)$ and $E_2=E(V_2)$ denote the corresponding strains. Then Richter's superposition postulate states that for $E=E(V_1V_2)$,
\begin{align}
	f(E_1)+f(E_2)=f(E)\label{eq:superpospri}
\end{align}
for some primary matrix function $f$, which depends on (and, in fact, determines) the specific choice of a strain mapping $F\mapsto E(F)$. This requirement is well known \cite{becker1893,agn_neff2014rediscovering,Hencky1928,Hencky1929,henckyTranslation,agn_neff2015geometry} to be satisfied for $f(\lambda)=\lambda$ and $E=\log V$, since%
\footnote{%
	It can easily be shown \cite{becker1893,agn_neff2014rediscovering} that under suitable normalization requirements, the \emph{only} strain tensor satisfying the condition $E(V_1V_2)=E(V_1)+E(V_2)$ for all coaxial stretches $V_1,V_2$ is the logarithmic Hencky strain $E(V)=E_0(V)=\log V$.%
}
\begin{align}
	\log(V_1V_2)=\log V_1+\log V_2\qquad\text{if}\quad V_1V_2=V_2V_1\,.\label{eq:knownstatement}
\end{align}
However, Richter's condition \eqref{eq:superpospri} is more general, allowing for an arbitrary choice of $f$. This generalization is what allows for any $E$ of the form \eqref{eq:straintensordefEulerian} to be considered a (Eulerian) strain tensor, since the representation
\begin{align}
	E=\ftilde(\log V)
\end{align}
implies that \eqref{eq:superpospri} is satisfied for $f=\ftilde\inv$. The somewhat unusual superposition principle \eqref{eq:superpospri} is thereby reduced to the better-known condition \eqref{eq:knownstatement}. As an example, consider again the Almansi strain tensor $E=\frac{1}{2}(\id-B\inv)$. Then for $\ftilde(\lambda)=\frac12(1-e^{-2\lambda})$ and $f(x)=\ftilde\inv(x)=-\frac12\log(1-2x)$,
\begin{alignat}{2}
	&\mathrlap{E=\frac{1}{2}(\id-B\inv)=\frac12 \big( \id-\exp(\log(V^{-2})) \big) = \ftilde(\log V)}\notag
	\\[.49em] \text{and}\qquad
	&f(E_1)+f(E_2) &&= f(\ftilde(\log V_1)) + f(\ftilde(\log V_2))
	\\ &&&= \log V_1 + \log V_2 = \log(V_1V_2) = f(\ftilde(\log(V_1V_2))) = f(E)
\end{alignat}
if $V_1V_2=V_2V_1$.

\subsection*{The strain deviator}
After giving a general definition of strain, Richter poses the following problem: given an arbitrary strain mapping $F\mapsto E(F)$, find an associated tensor valued mapping $F\mapsto D(F)$ that is invariant with respect to pure scaling transformations (i.e.\ $D(\lambda F)=D(F)$), reduces to $D=E$ if the deformation does not change the volume (i.e.\ $D(F)=E(F)$ if $\det F=1$) and coincides with the usual deviatoric strain tensor $\dev \varepsilon=\varepsilon-\frac{1}{3}\tr(\varepsilon)\cdot\id$ for infinitesimal deformations. From these conditions, Richter deduces the expression
\begin{equation}
\label{eq:deviatorDefiniton}
	D(F)=f\inv(\dev f(E(F)))=f\inv(\dev\log V)\,,
\end{equation}
where $f$ is given by \eqref{eq:superpospri} via the particular choice of the strain $E$. His deduction is based on the observation that the matrix logarithm naturally separates the isochoric and volumetric response, i.e.\ that
\begin{align}
	\log V=\dev(\log V)+\frac{1}{3}\tr(\log V)\cdot\id=\log\left(\frac{V}{\det V^{\afrac13}}\right)+\frac{1}{3}\log(\det V)\cdot\id\,.
\end{align}
In particular, if $D$ is defined by \eqref{eq:deviatorDefiniton}, then
\[
	D(\lambda F) = f\inv(\dev\log (\lambda V)) = f\inv(\dev\log V) = D(F)
\]
and, if $\det F = \det V = 1$,
\[
	D(F) = f\inv(\dev\log V) = f\inv(\log V) = E(F)\,.
\]

\subsection*{Richter's stress tensor}
In the following, we confine our attention to the setting of Cartesian coordinates. In that case, Richter proposes the use of the Cauchy stress tensor $\sigma$ and derives the necessary relations for the work corresponding to the displacement of surface elements. As a result, he obtains the formula
\begin{align}
	e^j\.\sigma&=\pdd{W}{j}\cdot\id+2\.\pdd{W}{k}\cdot L+3\.\pdd{W}{l}\cdot L^2\,, \label{Kirchhoff}
\end{align}
where $W(F)=W(j,k,l)$ is the isotropic energy potential in terms of the three invariants
\[
	j=\tr L\,, \quad k=\tr(L^2) \quad\text{and}\quad l=\tr(L^3)
\]
of the logarithmic strain $L=\log V$. Equation \eqref{Kirchhoff}, which had already been given by Richter in an earlier 1948 article \cite[page 207, eq.\ (3.9)]{richter1948}, can also be restated as a more common expression for the \emph{Kirchhoff stress} $\tau$ in hyperelasticity:
Using the notation
\[
	\widehat{W}(\log V) = W(F) = W(j,k,l) = W\big(\tr(\log V), \tr((\log V)^2),\tr((\log V)^3)\big)
\]
and the equalities
\begin{align*}
	D_{\log V}(j)&=D_L(\tr L)=\id\,,\\
	D_{\log V}(k)&=D_L(\tr(L^2))=D_L(\norm{L}^2)=2\.L\,,\\
	D_{\log V}(l)&=D_L(\tr(L^3))=3\.L^2\,,
\end{align*}
we find
\begin{align}
	D_{\log V}\widehat{W}(\log V)=\pdd{W}{j}\cdot\id+2\.\pdd{W}{k}\cdot L+3\.\pdd{W}{l}\cdot L^2\,.
\end{align}
Since
\begin{align}
	e^j=e^{\tr(\log V)}=e^{\log(\det V)}=\det V=\det F\,,
\end{align}
equation \eqref{Kirchhoff} can therefore be written as
\begin{equation}
\label{Kirchhoff2}
	\tau = \det F\cdot\sigma = D_{\log V}\widehat{W}(\log V)
	\,,
\end{equation}
where $\tau=\det F\cdot\sigma$ is the Kirchhoff stress tensor.
Formula \eqref{Kirchhoff} has been rediscovered several times \cite{Moreau76,vallee1978,hill1970,hill1978,ball2002openProblems} and is closely connected to Hill's inequality \cite{hill1970}, which is equivalent to
the condition that the elastic energy potential $W(F)=\widehat{W}(\log V)$ is convex with respect to the logarithmic strain tensor $\log V$. In particular, this convexity of $\widehat{W}$ is sufficient for $W$ to satisfy the Baker-Ericksen inequalities \cite{bakerEri54,buliga2002lower,silhavy1997mechanics}.

\bigskip

In the following, we provide a new translation of Richter's original 1949 article. For the sake of readability, the notation was updated to match more closely with current usage; a complete list of the changes made can be found in Table \ref{table:notation}.
The same updated notation has also been employed in translating the review of Richter's work by Ruth Moufang in \emph{Zentralblatt für Mathematik und ihre Grenzgebiete} as well as a Mathscinet review by William Prager. Apart from these notational changes, all equations as well as the equation numbering are identical to Richter's originally published version of the article. All numbered footnotes are part of the original article as well, whereas comments by the translators are marked as such.

\end{refsegment}
\begin{refsegment}

\numberwithin{equation}{section}
\newpage
\begin{center}
	{ \huge Strain tensor, strain deviator and stress tensor for finite deformations\\\medskip \large By \emph{Hans Richter} in Haltingen (Lörrach)\\[0.5cm]
	Zeitschrift für Angewandte Mathematik und Mechanik, Vol.~29, No.~3, 1949}
\end{center} 
\begin{abstract}
\par\medskip \noindent
	Postulates are laid down that have to be satisfied on forming the strain tensor, the strain deviator and the stress tensor, and thus the general form of these tensors are deduced in arbitrary coordinates. The mixed variant logarithmic strain tensor proves the simplest definition of the strain tensor. The deviator may be formed in the usual manner, and the invariants of it characterize the strain in an invariant way. If the stress tensor is defined accordingly, the form of the general law of elasticity continues to be invariant to coordinate transformations.
	\par\medskip \noindent
	Es werden Postulate aufgestellt, denen bei der Bildung des Verzerrungstensors, des Verzerrungsdeviators und des Spannungstensors zu gen\"{u}gen ist, und hieraus die allgemeine Gestalt dieser Tensoren in beliebigen Koordinaten abgeleitet. Als einfachste Definition des Verzerrungstensors erscheint die gemischt-variante logarithmische Deformationsmatrix, wo der Deviator in \"{u}blicher Weise gebildet werden kann, und wo die Invarianten des letzteren die Beanspruchung invariant charakterisieren. Bei entsprechender Definition des Spannungstensors bleibt die Gestalt des allgemeinen Elastizitätsgesetzes invariant gegen Koordinatentransformation.
	\par\medskip \noindent
	On \'{e}tablit des postulats pour la formation du tenseur de d\'{e}formation, du d\'{e}viateur de d\'{e}formation et du tenseur de tension. La forme g\'{e}n\'{e}rale de ces tenseurs en coordonn\'{e}es arbitraires en est d\'{e}duite. La matrice logarithmique (mixte-variante) de d\'{e}formation fournit la plus simple d\'{e}finition du tenseur de d\'{e}formation. Le d\'{e}viateur peut \^etre form\'{e} comme de coutume et ses invariantes caract\'{e}risent la sollicitation d'une mani\'{e}re invariante. Le tenseur de tension \'{e}tant d\'{e}fini conform\'{e}ment, la forme de la loi g\'{e}n\'{e}rale d'\'{e}lasticit\'{e} reste invariante dans toute transformation de coordonn\'{e}es.
\end{abstract}
\let\footnote\defaultfootnote%
\section{Introduction}

In the theory of finite elastic or plastic deformations, one generally considers the strain tensor which results from calculating the difference of the squares of the line elements in the deformed and initial state for general coordinates.\footnote{Cf.\ R.\ Moufang: \enquote{Volumtreue Verzerrungen bei endlichen Formänderungen}, Zeitschrift für Angewandte Mathematik und Mechanik 25/27 (1947), Pp.\ 209--214 \cite{moufang1947volumtreue}.} The use of this characterization of the state of strain is, of course, not compulsory. On the contrary, a more detailed analysis  shows that this usual definition of the strain tensor is not particularly well adapted to the problem of studying finite deformations. The problem of deducing a deviator, which only characterizes the change of shape without regarding the volume change, from the usual strain tensor already leads to peculiar difficulties and ambiguities \cite{moufang1947volumtreue}. The underlying reason for this is that the treatment of finite deformations has been approached too closely to the case of infinitesimal strains, where any deformation can be split into a pure stretch and a pure rotation by additive decomposition into a symmetric and a skew symmetric part. However, for finite deformations this additive decomposition is no longer possible; it is replaced by a multiplicative decomposition of the general deformation into a rotation and a stretch, with these factors no longer being commutative. Thus any attempt to establish definitions by additive decomposition must lead to fundamental difficulties.
\par\medskip\noindent
In this paper we want to proceed –- in a sense axiomatically –- by imposing on the necessary definitions certain a priori requirements we consider appropriate. Then, we demonstrate that among these admissible definitions, certain choices appear particularly natural.
\section{Notation and lemmas}
\subsection{Notation}
\begin{enumerate}
	\item 
	By Latin capital letters $A,B,\dotsc$ we denote elements of the space of $3\times3$-matrices.\footnote{Whether or not a matrix is a tensor is determined by (2.10).} $a_{ik}=(A)_{ik}$ is the entry in the $i$-th row and the $k$-th column. $\det A$ is the determinant of $A$. $\tr (A)$ is the trace of $A$, i.e.\ the sum of the elements on the main diagonal. $A^T$ is the matrix obtained by reflecting $A$ across its main diagonal. $\id$ is the identity matrix. $A\inv$ is the inverse of $A$.
	\item
	Latin lower case letters $x,y,\dotsc$ denote vectors: $x=(x_1,x_2,x_3)$. $\iprod{x,y}$ is the inner product. $x\times y$ is the cross product.
	\item
	$A\.x$ results from applying $A$ to $x$: $(A\.x)_i=\sum a_{ik}\.x_k$.
	\item
	Products $B\.A$ are read from right to left: %
	$(B\.A)\.x=B\.(A\.x)$.
	\item
	If $f(x)=\sum b_n \cdot x^n$, then, assuming convergence: $f(A)=\sum b_n\.A^n$;\\
	$\intd{f}(A)=f(A+\intd{A})-f(A)$, which coincides with $f'(A)\.\intd{A}$ only if $A\cdot\intd{A}=\intd{A}\cdot A$.\footnote{However, c.f. \ref{2.5}.}
	\end{enumerate}

\subsection{Lemmas}

\begin{enumerate}[label=(2.\arabic*), align=left, labelsep=0.25cm, labelwidth=\widthof{(2.10)},leftmargin=\labelwidth+\labelsep]
	\item
$\tr (A_1\.A_2 \cdots A_n)$ is invariant under cyclic permutations of the factors.
	\item
Each invariant of $A$ under affine transformation $A\to B\.A\.B\inv$ is a function of the three invariants $j=\tr(A)$, $k=\tr(A^2)$ and $l=\tr(A^3)$. The characteristic equation of $A$ is:
\begin{align*}
	\lambda^3-j\cdot \lambda^2+\half\.(j^2-k)\.\lambda-\left(\tel{3}\.l-\half\.j\.k+\tel{6}\.j^3\right)=0\,.
\end{align*}\label{2.2}
	\item
We have $f(B\.A\.B\inv)=B\.f(A)\.B\inv.$\label{2.3}
	\item
If $A$ has positive real eigenvalues, then $\log A$ is well defined and $\tr (\log A)=\log(\det A)$.\label{2.4}
	\item
If $B\.A = A\.B$, then $\tr (B\,\.\intd{f}(A))=\tr (B\.f'(A)\.\intd{A})$ even if $B\cdot\intd{A}=\intd{A}\cdot B$ does not hold.\label{2.5}
	\item
In Cartesian coordinates, a pure stretch $V$ is symmetric with positive eigenvalues.
	\item
In Cartesian coordinates, a Euclidean transformation $R$ satisfies $R\.R^T=\id$.	
	\item
Any $A$ with $\det A\neq 0$ can be uniquely represented in the form $A=V\cdot R$, i.e.\ as the composition mapping of a Euclidean transformation and a pure stretch. If $\det A>0$, then $R$ is a direct transformation, i.e.\ a pure Euclidean rotation. \label{2.8}
	\item
We have $\iprod{x, A\.y} = \iprod{y, A^Tx}.$	\label{2.9}
	\item
Let $y=\jacobi\.x$ be a coordinate transformation which maps $A$ onto $A^\#$. $A$ is a
		\begin{center}
			\begin{tabular}{ll}
twice-contravariant tensor if			&$A^\#=\jacobi\.A\.\jacobi^T\cdot (\det\jacobi)^n$,\\
twice-covariant tensor if					&$A^\#=(\jacobi\inv)^TA\.\jacobi\inv\cdot (\det\jacobi)^n$,\\
contravariant-covariant tensor if	&$A^\#=\jacobi\.A\.\jacobi\inv\cdot (\det\jacobi)^n$,\\
covariant-contravariant tensor if	&$A^\#=(\jacobi\inv)^TA\.\jacobi^T\cdot (\det\jacobi)^n$.\\
			\end{tabular}
		\end{center}

$A$ is called a proper tensor if $n=0$ holds; if $n\neq 0$, then $A$ is called a tensor density. (The coincidence of this somewhat uncommon representation of the tensor property with the usual one immediately results from symbolically setting $(A)_{ik}=x_i\.y_k$, where $x$ and $y$ are contravariant or covariant vectors). \label{2.10}
	\item
Let $x'=\jacobi\.x$ and $y'=\jacobi\.y$; then $x'\times y'=\det\jacobi\cdot (\jacobi\inv)^T(x\times y)$.\label{2.11}
\end{enumerate}
\section{The strain tensor}

We now consider which requirements can be imposed justifiably on the strain tensor. Afterwards we will study the feasibility of these requirements.
\par\medskip\noindent
Let $F$ be the matrix which maps the neighborhood of a point $\widehat{x}$ to the neighborhood of its image $x$ under $F$:
\begin{align}
	\intd{x}=F\.\intd{\widehat{x}}\,.
\end{align}
$F$ is the Jacobian matrix
\begin{align}
	(F)_{ik}=\pdd{x_i}{\widehat{x}_k}\,,\qquad\det F>0
\end{align}
and indicates the attained state of distortion. For plastic materials, where the state of stress does not only depend on the current state of distortion but also on the path leading to it, specifying only $F$ is not sufficient, whereas for elastic materials, $F$ suffices to characterize the distortion. For anisotropic materials the rotation contained in $F$ is essential as well. In this case, $F$ itself needs to be used for describing the strain, whereas every strain tensor which, like the common one, eliminates a Euclidean rotation is unsuitable. Consequently, such strain tensors are only meaningful for isotropic materials.

\subsection{Postulates}

Thus, under the explicit assumption of applicability to \emph{isotropic materials}, a strain tensor $E(F)$ associated with $F$ shall now be defined.\footnote{The notation $E(F)$ does not mean that $E$ is a function of $F$ in the sense of \ref{2.5}, but merely indicates that $E$ is associated with $F$.} Whereas $F$ is not a tensor since $F$ relates two different configurations, we want to require the tensor property for $E$. Hence, we obtain the first postulate:
\begin{postulat}
	\emph{$E$ is a tensor determined by $F$ and the matrices of the metric in $\widehat{x}$ and $x$.}
\end{postulat}
\noindent
Furthermore, the irrelevant rotation contained in $F$ shall be disregarded for $E$, i.e.\ $E$ shall not change if a Euclidean rotation $R$ is performed in $\widehat{x}$ \emph{prior} to the application of $F$. Instead, one could also require that a rotation being performed \emph{subsequent} to $F$ in $x$ shall not influence the strain tensor. This would imply that $F$ is considered a distortion in $\widehat{x}$ with a subsequent irrelevant rotation. We want to denote the tensor being associated with $\widehat{x}$ by $\widehat{E}$. The study of $E$ and $\widehat{E}$ is completely analogous and thus, in the following, we restrict ourselves to the study of $E$ and only mention the analogous results of $\widehat{E}$, where the corresponding quantities are marked by $\widehat{\hphantom{A}}$.
\par\medskip\noindent
The above property of $E$ and $\widehat{E}$ is expressed by the postulate
\begin{postulat}
	$E(F\.R)=E(F)$, \quad resp. \quad $\widehat{E}(R\.F)=\widehat{E}(F)$.
\end{postulat}
\noindent
Furthermore, we additionally require a \emph{superposition principle} for coaxial pure stretches via the postulate
\begin{postulat}
	\emph{Let $V_1$ and $V_2$ be two coaxial stretches: $V_1\,V_2=V_2\,V_1$. Let $E_1=E(V_1)$, $E_2=E(V_2)$ and $E=E(V_1\.V_2)$. Then there exists an invertible function $f(x)$ such that $f(E_1)+f(E_2)=f(E)$.} The function $f$ may depend on the coordinate system.
\end{postulat}
\noindent
Finally, we must require that the new definition transitions into the classical one for infinitesimal strains. This is ensured by the \emph{limit property}
\begin{postulat}
	\emph{For infinitesimal deformations $\id+\intd{F}$ in Cartesian coordinates the strain tensor turns into $\half\.(\intd{F}+(\intd{F})^T)+o(\intd{F})$}.\footnote{As usual, $y=o(x)$ means that: $\lim\frac{y}{x} =0\,.$}
\end{postulat}

\subsection{The realization of the postulates in Cartesian coordinates}

For the sake of simplicity, we first want to assume Cartesian coordinates. We denote an original point and its image under the deformation by $\opoint$ and $\ipoint$. The deformation matrix is now denoted by $\deform$. The corresponding strain tensors are $E_0$ and $\widehat{E_0}$.
\par\medskip\noindent
According to \ref{2.8} we can write
\begin{align}
	\deform=VR=R\.U\qquad\text{with}\qquad U=R\inv\.VR\,.
\end{align}
\noindent
To find this decomposition, we first consider the term $\deform\.\deformT$. For $x\neq 0$ we have: $0<\iprod{\deformT x,\deformT x}$, which, using \ref{2.9}, yields : $0<\iprod{x,\deform\.\deformT x}$. Thus, the symmetric matrix $\deform\.\deformT$ is positive definite and therefore obviously has a positive definite square root $V=\sqrt{\deform\.\deformT}$. Then $R$ can be represented in the form $R=V\inv\.\deform$. Correspondingly, we have $U^2=\deformT\deform$.
\par\medskip\noindent
By \textbf{V2}, we have $E_0(\deform)=E_0(V)$, resp. $\widehat{E_0}(\deform)=\widehat{E_0}(U)$. Therefore we can restrict ourselves to strain tensors which are defined for pure stretches.
\par\medskip\noindent
Now let $V$ be a pure infinitesimal stretch: $V=\id+\intd{V}$. Then by \textbf{V4} the equalities $E_0(\id+\intd{V})=\intd{V}+o(\intd{V})$ and $E_0(\id+\lambda\.\intd{V})=\lambda\.\intd{V}+o(\intd{V})$ hold for any positive number $\lambda$. Postulate \textbf{V3} then yields\\
$f(\intd{V}+o(\intd{V}))+f(\lambda\.\intd{V}+o(\intd{V}))=f((1+\lambda)\.\intd{V}+o(\intd{V}))$. Since this equation must hold for every $\lambda$ and $\intd{V}$, it follows that $f(x)=x+o(x)$ for $x$ sufficiently small\footnote{Since with $f$ every multiple of $f$ also satisfies postulate \textbf{V3}, $f'(0)$ can be normalized to $1$.}. Thus, if we set $Z=f(E_0)$, then for infinitesimal stretchings we obtain: $Z(\id+\intd{V})=\intd{V}+o(\intd{V})$.
\par\medskip\noindent
Now let $V$ again be a finite pure stretching. Then because of the positive eigenvalues of $V$ we can set:
\begin{align}
	L=\log V\,;\quad\text{resp.}\quad\widehat{L}=\log U\,:\quad\text{\glqq logarithmic strain tensor\grqq.}\label{3.4}
\end{align}
We then have: $\tel{n}\.L=\log\sqrt[n]{V}$ and thus for $n$ sufficiently large: $\sqrt[n]{V}=\id+\tel{n}\.L+o(\tel{n})$. Hence, $Z(\sqrt[n]{V})=\tel{n}\.L+o(\tel{n})$. In addition, we have $Z(V)=n\cdot Z(\sqrt[n]{V})=L+n\cdot o(\tel{n})$ by postulate \textbf{V3}. Since the left hand side of this equation is independent of $n$, we can let $n$ tend to infinity and obtain: $Z(V)=L$. In particular, this implies that $f(x)$ is uniquely determined up to an arbitrary factor.
\par\medskip\noindent
Consider the inverse function $f\inv$. Since $L$ is a uniquely invertible function of $V$ and consequently one of $V^2=\deform\.\deformT$, we finally have:
\begin{align}
	E_0=f\inv(L)=h(V)=k(\deform\.\deformT)\,;\quad\text{resp.}\quad\widehat{E_0}=f\inv(\widehat{L})=h(U)=k(\deformT\deform)\,.\label{3.5}
\end{align}
Conversely, the ansatz \eqref{3.5} always satisfies the postulates \textbf{V2} and \textbf{V3}, where $f$ is uniquely chosen as the inverse function of $f\inv$, whereas satisfying the limit condition \textbf{V4} requires that for small $x$ we have:
\begin{align*}
	f\inv(x)&=x+o(x)\,;\qquad h(1+x)=x+o(x)\,;\qquad k(1+x)=\half\.x+o(x)\,.\tag{3.5a}\label{3.5a}\\
\intertext{Indeed, we then have for infinitesimal deformations $\deform=\id+\intd{\deform}$:}
	\deform\.\deformT&=\id+(\intd{\deform}+\intd{\deform}^T)+o(\intd{\deform})\,,\quad V=\sqrt{\deform\.\deformT}=\id+\half\.(\intd{\deform}+\intd{\deform}^T)+o(\intd{\deform})\,,\\
	L&=\log V=\half\.(\intd{\deform}+\intd{\deform}^T)+o(\intd{\deform})\,.
\end{align*}
Thus for every $f\inv$ satisfying \eqref{3.5a}: $E_0=\half\.(\intd{\deform}+\intd{\deform}^T)+o(\intd{\deform})\,.$ 
\par\medskip\noindent
\emph{Every strain tensor being compatible with our postulates is thus identified with a function of the logarithmic strain tensor.} Based on our postulates, $E_0=L$ appears as the simplest definition of the strain tensor since here, the superposition principle is satisfied with $f(x)\equiv x$. As we shall later see, this definition will also appear as the simplest one for taking the deviator.
\par\medskip\noindent
If, in addition, a Euclidean rotation $R_1$ is performed subsequently to $F$, then because of $R_1\.F=R_1\.V\.R=R_1\.V\.R_1\inv\cdot R_1\.R$ the stretch $V$ turns into $R_1\.V\.R_1\inv$. We obtain the same transition if a Euclidean coordinate transformation $\ipoint_1=R_1\.\ipoint$ is performed. According to \ref{2.3}, $E_0$ then turns into $h(R_1\.V\.R_1\inv)=R_1\.E_0\.R_1\inv$. Thus the axes of $E_0$ are simply rotated along for subsequent application of $R_1$. If we identify the last formula with the result of a coordinate transformation, we conclude from \ref{2.10} that $E_0$ transforms like a tensor; since $R_1\inv=R_1^T$, there is no distinction with respect to co-contra-variance. Clearly, we obtain a corresponding result for $\widehat{E_0}$.

\subsection{Extension to curvilinear coordinates}

We now proceed from Cartesian coordinates $\ipoint$ to arbitrary coordinates $x$: $x=x(\ipoint)$. For a neighborhood of the undeformed material let $\intd{\widehat{x}}=\widehat{\jacobi}\.\intd{\widehat{\ipoint}}$, for the corresponding neighborhood in the deformed material let $\intd{x}=\jacobi\.\intd{\ipoint}$. $\widehat{\jacobi}$ and $\jacobi$ are the Jacobian matrices of $x=x(\ipoint)$ in $\widehat{x}$ and $x$, respectively.
\par\medskip\noindent
For a line element in $x$ we obtain, using \ref{2.9}: $\intd{\ipoint}=\iprod{\jacobi\inv\intd{x},\jacobi\inv\intd{x}}=\iprod{\intd{x},(\jacobi\inv)^T\jacobi\inv\intd{x}}$. Hence,
\begin{align*}
	G&=G^T=(\jacobi\.\jacobi^T)\inv\tag{3.6a}\\
\intertext{is the matrix of the metric in $x$. Correspondingly,}
	\widehat{G}&=(\widehat{\jacobi}\.\widehat{\jacobi}^T)\inv\tag{3.6b}
\end{align*}
defines the metric in $\widehat{x}$.
\par\medskip\noindent
The deformation of the material now appears as: $\intd{x}=\jacobi\.\deform\.\widehat{\jacobi}\inv\intd{\widehat{x}}$. Thus, $\deform$ is changed to
\stepcounter{equation}
\begin{align}
	F=\jacobi\.\deform\.\widehat{\jacobi}\inv\,,\qquad (F)_{ik}=\pdd{x_i}{\widehat{x}_k}\,.
\end{align}
Conversely,
\begin{align*}
	\deform=\jacobi\inv F\.\widehat{\jacobi}\qquad\text{and}\qquad\deformT=\widehat{\jacobi}^TF^T(\jacobi\inv)^T\,,\tag{3.7*}\label{3.7*}
\end{align*}
from which we immediately obtain:
\begin{align}
	\left\{\begin{alignedat}{2}
		\quad V^2&=\deform\.\deformT&&=\jacobi\inv F\.\widehat{G}\inv F^T(\jacobi\inv)^T\\
		U^2&=\deformT\deform&&=\widehat{\jacobi}^TF^TG\.F\.\widehat{\jacobi}\,.
	\end{alignedat}\right.\label{3.8}
\end{align}
Using the last two formulae, the matrices $\deform$, $V$ and $U$ associated with Cartesian coordinates can be expressed in terms of $F$ and the transformation matrices $\jacobi$ and $\widehat{\jacobi}$.

\subsubsection{Case of the non-mixed tensor}

We first assume that the strain tensor $E$ is defined twice-contravariant and satisfies the postulates \textbf{V1}-\textbf{V4}. Then \ref{2.10} implies:
\begin{align*}
	E=\jacobi\.E_0\.\jacobi^T\,,
\end{align*}
where $E_0$ is one of the tensors from \eqref{3.5}.
\par\medskip\noindent
To study the particular shape of $E_0$, we consider the special case where $\deform$ is a pure stretch $V$ in the coordinate axes and coaxial to $\jacobi$. Hence
\begin{center}
	$V=\matr{\lambda_1&&0\\&\lambda_2&\\0&&\lambda_3}$ \qquad and \qquad $\jacobi=\matr{\varrho_1&&0\\&\varrho_2&\\0&&\varrho_3} .$
\end{center}
Then because of \eqref{3.5} $E$ is again given in principal axis and has the eigenvalues $\varrho_v\cdot h(\lambda_v)$. The superposition principle now requires the existence of a function $f(x)$, whose coefficients may contain the $\varrho_v$, such that:
\begin{align}
	f(\varrho_v^2\.h(\lambda_v))+f(\varrho_v^2\.h(\mu_v))=f(\varrho_v^2\.h(\lambda_v\.\mu_v))	
\end{align}
for arbitrary $\lambda_v$ and $\mu_v$. Therefore,
\begin{align*}
	f(\varrho_v^2\.h(\lambda))=C_v\cdot\log\lambda\,,\qquad C_v=C_v(\varrho_1,\varrho_2,\varrho_3)\,.
\end{align*}
By differentiation with respect to $\lambda$ we obtain 
\begin{align}
	\varrho_v^2\.f'(\varrho_v^2\.h(\lambda))\cdot h'(\lambda)=C_v\cdot\tel{\lambda}\,.
\end{align}
In particular, if we set $\lambda=1$, then \eqref{3.5a} implies $\varrho_v^2\cdot f'(0)=C_v$. Therefore the normalization  $f'(0)=1$ yields:
\begin{align*}
	f'(\varrho_v^2\.h(\lambda))=\tel{\lambda\cdot h'(\lambda)}\,.
\end{align*}
The right-hand side of this equation is independent of $\varrho_v$. We must therefore have $f'(x)\equiv f'(0)=1$, which implies $h(\lambda)=\log\lambda$. Then $E_0=\log V=L$ and thus
\begin{align*}
	E=\jacobi\.L\.\jacobi^T .
\end{align*}
Conversely, this definition of $E$ satisfies all postulates \textbf{V1}-\textbf{V4} for arbitrary $\deform$ and $\jacobi$, since the superposition principle is purely additive for $L$ and therefore transfers to an additive law in terms of $E$ for multiplication with $\jacobi$ from the left and with $\jacobi^T$ from the right.
\par\medskip\noindent
\emph{Hence, there is only one possibility to define a non-mixed tensor $E$ such that our postulates are satisfied, namely: $E=\jacobi\.L\.\jacobi^T$.}
\par\medskip\noindent
Since $L=\log V=\half\.\log V^2$, and due to \eqref{3.8}, we finally obtain
\begin{align}
	E&=\half\.\jacobi\.\log(\jacobi\inv F\.G\inv F\.\jacobi\inv)\.\jacobi\,.\\
\intertext{Correspondingly, we find}
	\widehat{E}&=\half\.\widehat{\jacobi}\log(\widehat{\jacobi}\.F^TG\.F\.\widehat{\jacobi})\.\widehat{\jacobi}\,.\notag
\end{align}
If we expand the logarithm for sufficiently small stretches into a power series, then we will see that $E$ and $\widehat{E}$ indeed only depend on $F$, $\widehat{G}$ and $G$. However, the representation by these matrices is very inconvenient. Moreover, the invariants of $E$ are different from those of $E_0$. This is unpleasant because e.g.\ in the theory of elasticity of finite deformations it must be assumed that the thermodynamic quantities like internal energy, entropy etc.\ are functions of the invariants of strain. Now, if these quantities are changed under coordinate transformations, additional difficulties will emerge. Then it is also no longer possible to describe the character of the deformation independently of the choice of coordinates by using the invariants (cf.\ chapter 4).
\par\medskip\noindent
Of course, the corresponding considerations hold also for the case where $E$ or $\widehat{E}$ is twice-covariant. Therefore it will be sufficient to waive the symmetry advantage being associated with non-mixed tensors.

\subsubsection{Case of the mixed tensor}

In the case where $E$ is covariant-contravariant, \ref{2.10} implies: $E=(\jacobi\inv)^TE_0\.\jacobi^T$. Because of \ref{2.3}, $E$ automatically satisfies the superposition principle with the same $f(x)$ as $E_0$. In particular, the uniquely determined, normalized $f(x)$ is independent of the choice of coordinates. Furthermore, $E$ has the same invariants as $E_0$. Every function of $E$, whose coefficients depend on the invariants of $E$, transforms to the same function of $E_0$.
\par\medskip\noindent
From the simplest definition $E_0=L$ we now obtain for arbitrary coordinates: $L^{\ast}=(\jacobi\inv)^TL\.\jacobi^T$ or, because of \eqref{3.4} and \eqref{3.8}:
\begin{align}
	L^{\ast}&=\half\.(\jacobi\inv)^T\log(\jacobi\inv F^T\widehat{G}\inv F^T(\jacobi\inv)^T)\.\jacobi^T\notag\\
\intertext{or}
	L^{\ast}&=\half\.\log(G\.F\.\widehat{G}\inv F^T)\,:\quad\text{\glqq logarithmic strain tensor\grqq.}\label{3.12}
\end{align}
$L^{\ast}$ satisfies the superposition principle with $f(x)\equiv x$.
\par\medskip\noindent
\emph{The most general strain tensor satisfying our postulates is then given by
\begin{align}
	E=f\inv(L^{\ast})=h\!\left(\sqrt{G\.F\.\widehat{G}\inv F^T}\right)=k(G\.F\.\widehat{G}\inv F^T)\,,\label{3.13}
\end{align}
where $f\inv$, $h$ and $k$ satisfy the conditions \eqref{3.5a}.}
\par\medskip\noindent
Completely analogous, one obtains
\begin{align*}
	\widehat{E}=f\inv(\widehat{L}^{\ast})=h\!\left(\sqrt{F^TG\.F\.\widehat{G}\inv}\right)=k(F^T\.G\.F\.\widehat{G}\inv)\,,\quad\widehat{L}^{\ast}=\half\.\log(F^TG\.F\.\widehat{G}\inv)\,.
\end{align*}
Up to an arbitrary factor, the function $f(x)$ of the superposition principle is the inverse function of $x=f\inv(y)$.
\par\medskip\noindent
In the case where $E$ and $\widehat{E}$ are contravariant-covariant, it is convenient to proceed correspondingly:
\begin{align*}
	\begin{cases}
		&L^{\ast}=\half\.\log(F\.\widehat{G}\inv F^TG)\\
		\quad\text{and}\qquad&\tag{3.12a}\\
		&{\widehat{L}}^{\ast}=\half\.\log(\widehat{G}\inv F^TG\.F)\,.
	\end{cases}\label{3.12a}
\end{align*}
Every other relation remains unchanged.

\subsection{Computation of the dilatation $v$}

The dilatation being associated with $F$ is $v=\det \deform$; thus, with \eqref{3.7*}:
\begin{align*}
	v=(\det\jacobi)\inv\cdot\det{\widehat{\jacobi}}\cdot\det{F}\,.	
\end{align*}
However, (3.6), \eqref{3.12} and \eqref{3.12a} yield:
\begin{align*}
	\det (e^{2\.L^*})=\det G\cdot (\det\widehat{G})\inv\cdot (\det F)^2=v^2\,.	
\end{align*}
Hence, due to \ref{2.4}:
\stepcounter{equation}
\begin{align*}
	\log v&=\tr(L^*)=\tr(\widehat{L}^*)\tag{3.14a}\label{3.14a}\\
\intertext{or by \eqref{3.13}}
	\log v&=\tr(f(E))=\tr(f(\widehat{E}))\,.\tag{3.14b}
\end{align*}
	
\subsection{Relation to the usual strain tensor}

The usual definition\transcomment{Here, $T=\frac{1}{2}(\id-B\inv)$, $\widehat{T}=\frac{1}{2}(C-\id)$.} of the strain tensor $T$, resp.\ $\widehat{T}$, is\footnote{See e.g.\ \cite{moufang1947volumtreue}}
\begin{align*}
	\intd{s}^2-\intd{\widehat{s}}^2=2\iprod{\intd{x}, T\.\intd{x}}=2\.\iprod{\intd{\widehat{x}},\widehat{T}\.\intd{\widehat{x}}}\,.
\end{align*}
Now, together with \ref{2.9} we get
\begin{align*}
	\intd{s}^2&=\iprod{\intd{x},G\.\intd{x}}=\iprod{F\.\intd{\widehat{x}},G\.F\.\intd{\widehat{x}}}=\iprod{\intd{\widehat{x}},F^TG\.F\.\intd{\widehat{x}}}\\
\intertext{and correspondingly}
	\intd{\widehat{s}}^2&=\iprod{\intd{x},(F\inv)^T\widehat{G}\.F\inv\intd{x}}=\iprod{\intd{\widehat{x}},\widehat{G}\.\intd{\widehat{x}}}\,.
\end{align*}
Hence
\begin{align*}
	2\.T=G-(F\.\widehat{G}\inv F^T)\inv\qquad\text{and}\qquad 2\.\widehat{T}=F^TG\.F-\widehat{G}\,.	
\end{align*}
In order to identify the type of co-contra-variance, we use \eqref{3.8} to rewrite:
\begin{align}
	2\.T&=(\jacobi\inv)^T\cdot (\id-\jacobi^T(F\inv)^T\widehat{G}\.F\inv\jacobi)\cdot\jacobi\inv=(\jacobi\inv)^T\cdot (\id-V^{-2})\cdot\jacobi\inv\label{3.15}\\
\intertext{and correspondingly}
	2\.\widehat{T}&=(\widehat{\jacobi}\inv)^T\cdot (U^2-\id)\.\widehat{\jacobi}\inv.\tag{3.15a}
\end{align}
Thus, according to \ref{2.10}, $T$ and $\widehat{T}$ are twice-covariant symmetric tensors. The superposition principle is not satisfied for these tensors and the invariants change under coordinate transformation. However, the combined tensors $T\.G\inv$, $G\inv T$, $\widehat{G}\inv\widehat{T}$ and $\widehat{T}\.\widehat{G}\inv$ satisfy all the established postulates \textbf{V1}-\textbf{V4}. From \eqref{3.15} we infer for the superposition principle that one has to set\transcomment{$f(T)=f(\frac{1}{2}(\id-B\inv)=-\frac{1}{2}\log(\id-2\cdot\frac{1}{2}(\id-B\inv))=-\frac{1}{2}\log(B\inv)=\log V$ for $f(x)=-\frac{1}{2}\log(1-2x)$.} $f(x)=-\half\.\log(1-2x)$ with respect to $T\.G\inv$ and $G\inv T$, but $f(x)=\half\.\log(1+2x)$ with respect to $\widehat{G}\inv\widehat{T}$ and $\widehat{T}\.\widehat{G}\inv$.
\par\medskip\noindent
Hence, with (3.14),
\begin{align*}
	v^2&=(\det(\id-2\.T\.G\inv))\inv=(\det(\id-2\.G\inv T))\inv\\
	&=\det(\id+2\.\widehat{G}\inv\widehat{T})=\det(\id+2\.\widehat{T}\.\widehat{G}\inv)\,
\end{align*}
for the dilatation $v$.
\section{The strain deviator}
\subsection{Postulates}

The strain deviator $D$ shall be derived from the strain tensor and only characterize the change of shape associated with the deformation. The required postulates immediately follow:
\begin{postulat2}
	\emph{If two deformations differ only by a scaling, then they have the same $D$.}
\end{postulat2}
\begin{postulat2}
	\emph{If the deformation does not change the volume, then $D=E$.}
\end{postulat2}

\subsection{Realization of the postulates}

A scaling in the undeformed or deformed state has the form $\lambda\.\id$, $\lambda>0$, with the volume dilatation $\lambda^3$. If we set
\begin{align*}
	F=v^\tel{3}\.\id\cdot v^{-\tel{3}}F=v^\tel{3}\.\id\cdot F_1\,,
\end{align*}
then postulate \textbf{D1} yields: $D(F)=D(F_1)$. Since $F_1$ is not associated with a volume dilatation, we have: $D(F)=E(F_1)=f\inv(L_1^{\ast})$, where $L_1^{\ast}=-\tel{3}\.\log v\cdot\id+L^\ast$ according to \eqref{3.12} and \eqref{3.12a}. Using \eqref{3.13} and \eqref{3.14a} we conclude that
\begin{align*}
	D=f\inv(L^\ast-\tel{3}\.\tr L^\ast\cdot\id)=f\inv(f(E)-\tel{3}\.\tr f(E)\cdot\id)\,.	
\end{align*}
The common deviator of a matrix $A$ is denoted by
\begin{align}
	\dev A=A-\tel{3}\.\tr A\cdot\id\,.
\end{align}
With this notation, we can reformulate the strain deviator as:
\begin{alignat}{2}
	D&=f\inv(\dev L^\ast)&&=f\inv(\dev f(E))\,.\label{4.2}\\
\intertext{Correspondingly,}
	\widehat{D}&=f\inv(\dev\widehat{L}^{\ast})&&=f\inv(\dev f(\widehat{E}))\,.\notag
\end{alignat}
Conversely, the postulates \textbf{D1} and \textbf{D2} are obviously satisfied for this definition as well. If $F$ is multiplied by $\lambda >0$, then $L^\ast$ turns into $L^\ast+\log\lambda\cdot\id$. Thus $\dev L^\ast$ and consequently $D$ are left unchanged. If additionally $v=1$, then \eqref{3.14a} implies $\tr L^\ast=0$, hence $L^\ast=\dev L^\ast$ and consequently $D=f\inv(L^\ast)=E$. Note also that $D$ is automatically a tensor if we use $E$ as a mixed tensor. This observation suggests a preference towards mixed tensors over non-mixed tensors.
\par\medskip\noindent
Taking the deviator is simplest for $E=L^\ast$, where $D=L^\ast$. Thus the use of the logarithmic strain tensor also allows the common deviator procedure for arbitrary coordinates.
\par\medskip\noindent
It should additionally be noted that for infinitesimal strains in Cartesian coordinates the new notion of the deviator turns into the original one. If $\deform=\id+\intd{\deform}$ is an infinitesimal deformation, then\\ $L=\half\.(\intd{\deform}+(\intd{\deform})^T)+o(\intd{\deform})$, therefore \eqref{4.2}, together with \eqref{3.5a}, yields
\begin{align*}
	D=\dev L+o(\dev L)=\dev\left(\half\.\left(\intd{\deform}+(\intd{\deform})^T\right)\right)+o(\intd{\deform})\,.
\end{align*}
For the common mixed strain tensor $T\.G\inv$ we found in chapter 3.5 that $f(x)=-\half\.\log(1-2\.x)$, thus $f\inv(y)=\half\.(1-e^{-2\.y})$ and $v=(\det(\id-2\.T\.G\inv))^{-\half}$. Moreover, $L^\ast=-\half\.\log(\id-2\.T\.G\inv)$ and hence $2\.\dev L^\ast=-\log(\id-2\.T\.G\inv)-\frac{2}{3}\.\log(v)\cdot\id$. Therefore, we finally obtain:\footnote{Cf.\ the somewhat different construction by Moufang \cite{moufang1947volumtreue}.}
\begin{align}
	D&=(\det (\id-2\.T\.G\inv))^{-\tel{3}}\cdot\left(T\.G\inv-\half\cdot\left[1-\sqrt[3]{\det (\id-2\.T\.G\inv)}\right]\cdot\id\right)\,.\\
\intertext{Correspondingly,}
	\widehat{D}&=(\det(\id+2\.\widehat{T}\.\widehat{G}\inv))^{-\tel{3}}\cdot\left(\widehat{T}\.\widehat{G}\inv-\half\cdot\left[\sqrt[3]{\det(\id+2\.\widehat{T}\.\widehat{G}\inv)}-1\right]\cdot\id\right)\,.\tag{4.3a}
\end{align}

\subsection{The strain invariants}

To characterize the state of strain through invariants we choose the dilatation (or a function of the same) as the first suitable invariant of $E$, whereas the other two invariants characterize the change of shape, i.e.\ they shall be left unchanged by additional scaling. Since for the use of the mixed tensors --- which is assumed in the following --- every invariant of $E$ is also an invariant of $L^\ast$, we can choose $\tr L^\ast$ as the first invariant by \eqref{3.14a}. According to section 2, the other two invariants must be invariants of $\dev L^\ast$. From this we conclude that the state of strain is characterized by
\begin{align}
	\left\{\begin{alignedat}{2}
		 j&=\tr L^\ast&&\text{for the dilatation}\\
		 y&=\tr((\dev L^\ast)^2)\,,\quad z=\tr ((\dev L^\ast)^3)\qquad&&\text{for the change of shape.}
	\end{alignedat}\right.\label{4.4}
\end{align}
Since $L^\ast=(U\inv)^TL\.U^T$, we have $y=\tr((\dev L)^2)$ and $z=\tr((\dev L)^3)$, therefore $y$ and $z$ characterize the change of shape independently of the choice of coordinates.
\par\medskip\noindent
Because of $\tr(\dev L)=0$, the characteristic equation of $\dev L$ according to \ref{2.2} is
\begin{align}
	x^3-\frac{y}{2}\.x-\frac{z}{3}=0\,.\label{4.5}
\end{align}
In order for this equation to have three real roots, the quantity
\begin{align}
	\zeta=\frac{z^2}{y^3}
\end{align}
must satisfy the condition
\begin{align}
	0\leq\zeta\leq\tel{6}.
\end{align}
The geometrical meaning of $\zeta$ results from the following observation. Let $V$ be an arbitrary pure stretch. Then we can identify $V$ with the $n$-fold application of the pure stretch $V_n=\sqrt[n]{V}$. Here, $L_n=\log\sqrt[n]{V}=\tel{n}\cdot\log V=\tel{n}\.L$; thus $\dev L_n=\tel{n}\.\dev L$ and consequently $y_n=\tel{n^2}\cdot y$ and $z_n=\tel{n^3}\.z$. From this we infer that $\zeta_n=\frac{z_n^2}{y_n^3}=\zeta$. --- Conversely, if $\zeta_1=\zeta_2$ for two stretches $V_1$ and $V_2$, then $y_1=\lambda^2\.y$ and $z_1=\lambda^3\.z$. Then according to \eqref{4.5}, the eigenvalues of $V_1$ are the $\lambda$-th power of the eigenvalues of $V_2$. Thus, disregarding a possible rotation, we have $V_1=V_2^{\lambda}$. Hence we can think of $V_1$ and $V_2$, up to a modification of the principal axes, as resulting from the same infinitesimal stretch (using the inverse for negative $\lambda$). This means that $\zeta$ determines the \emph{character} of the deformation.
\par\medskip\noindent
The uniaxial and volume preserving stretch is represented in suitably rotated Cartesian coordinates by
\begin{align*}
	V=\matr{\lambda&0&0\\0&\lambda^{-\half}&0\\0&0&\lambda^{-\half}} .
\end{align*}
Then \quad\qquad\qquad\qquad\qquad\qquad $L=\dev L=\log\lambda\cdot\matr{1&0&0\\0&-\half&0\\0&0&-\half}$.\\
Thus $y=\log^2\lambda\cdot\frac{3}{2}$ and $z=\log^3\lambda\cdot\frac{3}{4}$, which results in $\zeta=\tel{6}$.\\
On the other hand, we obtain for a volume preserving simple shear
\begin{align*}
	\deform=\matr{1&\lambda&0\\0&1&0\\0&0&1}\qquad\text{and thus}\qquad V^2=\deform\.\deformT=\matr{1+\lambda^2&\lambda&0\\\lambda&1&0\\0&0&1} .
\end{align*}
For the eigenvalues of $V^2$, the characteristic equation yields: $\lambda_1=1$, $\lambda_2\cdot\lambda_3 =1$.
For the eigenvalues of $L$ we thus have: $\mu_1=0$, $\mu_2+\mu_3=0$. This implies $y>0$, $z=0$; therefore $\zeta =0$. Hence, we have found that:
\par\medskip\noindent
\emph{$\zeta =\frac{z^2}{y^3}$ determines the character of the deformation. The extreme value $\zeta =0$ corresponds to simple shearing and the other extreme value $\zeta=\tel{6}$ to uniaxial stretching.}
\par\medskip\noindent
The \emph{amount} of change of shape at infinitesimal strains is usually characterized by $\sqrt{\tr D^2}$. We have just shown that $D\approx\dev L$ for infinitesimal deformations, hence $\sqrt{y}$ is identified with the amount of change of shape at infinitesimal deformations.
\par\medskip\noindent
On the other hand, if $V$ is a finite scaling, then $\sqrt{y}=n\cdot\sqrt{y_n}$ as demonstrated above. Since for sufficiently large $n$, $y_n$ represents the amount of change of shape for the infinitesimal strain $\sqrt[n]{V}$, it is reasonable to use $\sqrt{y}=n\.\sqrt{y_n}$ as a measure for the amount of change of shape resulting from an $n$-fold application of $\sqrt[n]{V}$, i.e.\ for $V$. From this we finally conclude: 
\par\medskip\noindent
\emph{$\sqrt{y}$ characterizes the amount of change of shape.}
\section{The stress tensor}

The stress tensor $\widetilde{\sigma}$ must characterize the state of stress in the point $x$ of the deformed configuration such that for a suitable definition of a surface element $\intd{A}$ in $x$, the force acting on $\intd{A}$ is given by $\widetilde{\sigma}\.\intd{A}$. Even though the components of $\widetilde{\sigma}$ can, of course, be expressed in the coordinates of $\widehat{x}$ as well by using the transformation formulae (transition into Lagrangian coordinates), $\widetilde{\sigma}$ remains associated with $\intd{A}$. The attempt to directly connect the stresses with $\intd{\widehat{A}}$ in the reference configuration, i.e.\ to construct a $\widehat{\sigma}$, is unnatural from a physical point of view.\transcomment{The stress $\widehat{\sigma}$ described here, connecting the surface element $\intd{\widehat{A}}$ to the occurring forces, corresponds to the first Piola-Kirchhoff stress.} We will therefore refrain from such an approach.

\subsection{Postulates}

For Cartesian coordinates, the stress matrix $\sigma$ yields the force $\intd{\force_0}$ acting on a surface element $\intd{A_0}$ in the point $\ipoint$ in the form: $\intd{\force_0}=\sigma\.\intd{A_0}$. In general, it can be assumed that external forces acting on the material do not generate volume dependent torques. Then it is well known that $\sigma$ is symmetric. However, this symmetry does not need to be assumed in the following.
\par\medskip\noindent
For arbitrary curvilinear coordinates, it is necessary to define a surface element $\intd{A}$ suitably as the transformed element of $\intd{A_0}$. Then the stress tensor $\widetilde{\sigma}$ shall be constructed such that the force acting on the surface element is again given by $\widetilde{\sigma}\.\intd{A}$. Applying a translation by the vector $\intd{z}$ to the surface element corresponds to the work $\iprod{\intd{z},\widetilde{\sigma}\.\intd{A}}$. From this we deduce the following postulates:
\begin{postulat3}
	\emph{$\widetilde{\sigma}$ is a tensor or a tensor density.}
\end{postulat3}
\begin{postulat3}
	\emph{For the surface element we have $\intd{A}=\const\.\intd{A}_0$, where $\const$ must be chosen suitable.}
\end{postulat3}
\begin{postulat3}
	\emph{If the surface element is displaced by $\intd{z}$, then the corresponding work is $\intd{\work}=\iprod{\intd{z},\widetilde{\sigma}\.\intd{A}}$.}
\end{postulat3}
	
\subsection{The realization of the postulates}

As a numerical quantity, $\intd{\work}$ must be invariant under coordinate transformation. Hence
\begin{align}
	\iprod{\intd{z},\widetilde{\sigma}\.\intd{A}}=\iprod{\intd{z}_0,\sigma\.\intd{A_0}}\label{5.1}
\end{align}
if $\intd{z_0}=\jacobi\inv\intd{z}$ is the corresponding translation vector in Cartesian coordinates. Now we obtain with postulate \textbf{P2} and \ref{2.9}:
\begin{align*}
	\iprod{\intd{z_0},\sigma\.\intd{A_0}}=\iprod{\jacobi\inv\intd{z},\sigma\.\const\inv\.\intd{A}}=\iprod{\intd{A},(\const\inv)^T \sigma^T\.\jacobi\inv\intd{z}}=\iprod{\intd{z}, (\jacobi\inv)^T \sigma\.\const\inv\intd{A}}.
\end{align*}
Since $\intd{z}$ and $\intd{A}$ are arbitrary vectors, the comparison with \eqref{5.1} yields:
\begin{align}
	\widetilde{\sigma}=(\jacobi\inv)^T\sigma\.\const\inv.\label{5.2}
\end{align}
Equation \ref{2.10} indicates that $\widetilde{\sigma}$ is either \quad $(\alpha)$ twice-covariant for $\const=\jacobi\cdot\sqrt{(\det G)^n}$ \quad or \quad $(\beta)$ covariant-contravariant for $\const=(\jacobi\inv)^T\cdot\sqrt{(\det G)^n}$.
\par\medskip\noindent
In fact, $\intd{A}$ is contravariant in case $(\alpha)$ and covariant in case $(\beta)$. If we choose the length of $\intd{A}$ as the geometrical quantity of the surface element, then we have to set $n=0$. As a consequence, $\widetilde{\sigma}$ is a proper tensor. Namely, in the cases $(\alpha)$ and $(\beta)$ we have:
\begin{align*}
	\widetilde{\sigma}&=(\jacobi\inv)^T\sigma\.\jacobi\inv\tag{5.2$\alpha$}\\
\text{and}\qquad\qquad&\\
	\widetilde{\sigma}&=(\jacobi\inv)^T\sigma\.\jacobi^T.\tag{5.2$\beta$}
\end{align*}
In case $(\alpha)$, we then have $\intd{A}=\jacobi\.\intd{A_0}$. If the surface element $\intd{A_0}$ is generated by the vectors $\intd{x_{10}}$ and $\intd{x_{20}}$, then $\intd{A_0}=\intd{x_{10}}\times\intd{x_{20}}=\jacobi\inv\intd{x_1}\times \jacobi\inv\intd{x_2}$ or using \ref{2.11}: $\intd{A_0}=\sqrt{\det G}\cdot \jacobi^T(\intd{x_1}\times\intd{x_2})$. Hence
\begin{align*}
	\intd{A}&=\sqrt{\det G}\cdot G\inv(\intd{x_1}\times\intd{x_2})\,.\tag{5.3$\alpha$}\\
\intertext{On the other hand, in case $(\beta)$ we obtain}
	\intd{A}&=\sqrt{\det G}\cdot(\intd{x_1}\times\intd{x_2})\,.\tag{5.3$\beta$}
\end{align*}
Clearly, it does not matter whether one prefers to use contravariant or covariant $\intd{A}$ for calculations. As shown by (5.3), however, the covariant definition $(\beta)$ yields the simpler formula, although in this case, symmetry of $\widetilde{\sigma}$ does not follow from the symmetry of $\sigma$. On the other hand, all invariants of $\widetilde{\sigma}$ still remain unchanged under coordinate transformation.

\subsection{The power for infinitesimal strains}

Now we assume that in a spatial neighborhood of $\ipoint$, a homogeneous state of stress defined by $\sigma$ occurs. Suppose that a closed volume $V$ has the boundary surface $\mathcal{F}$ with the surface elements $\intd{A_0}$. We now apply a homogeneous infinitesimal deformation $\id+\intd{F}$ in the neighborhood of $\ipoint$. As a result, the surface element $\intd{A_0}$ is displaced by the vector $\intd{F}\.r_0$, provided that $r_0$ was its original distance to the origin. Since, due to symmetry, the simultaneous infinitesimal rotation and distortion of the surface element do not require any power, the entire work with respect to the volume is given by
\begin{align*}
	V\cdot\intd{\work}&=\iint\limits_V\iprod{\intd{F}\.r_0,\sigma\.\intd{A}_0}=\iint\limits_V\iprod{\sigma^T\.\intd{F}\.r_0,\intd{A_0}}\\
	&=\iiint\limits_V\text{div} (\sigma^T\.\intd{F}\.r_0)\.\intd{V}=\iiint\limits_V\tr(\sigma^T\.\intd{F})\.\intd{V}\,.
\end{align*}
Hence, the work per unit volume is
\begin{align}
	\intd{\work}=\tr(\sigma^T\.\intd{F})\tag{5.4}\,.\label{5.4}
\end{align}
If $\id+\intd{F}$ is an infinitesimal radial scaling, i.e.\ $\intd{F}=\intd{\lambda}\cdot\id$ with the dilatation $\frac{\intd{V}}{V}=3\cdot\intd{\lambda}$, then $\intd{\work}=\intd{\lambda}\cdot\tr \sigma^T=\frac{\intd{V}}{V}\cdot\tel{3}\.\tr \sigma^T$. If the hydrostatic stress $\hpressure$ occurs, then $\intd{\work}=\frac{\intd{V}}{V}\cdot\hpressure$. For non-symmetric $\sigma$, the hydrostatic stress is represented by $\tel{3}\tr \sigma^T$ as well, which is why $\tel{3}\.\tr \sigma^T$ is called the mean stress $\hpressure$.
\par\medskip\noindent
For arbitrary coordinates, according to \eqref{5.2}, the mean stress is given by
\begin{align*}
	\hpressure=\tel{3}\.\tr \sigma^T=\tel{3}\.\tr(\const^T\widetilde{\sigma}^T\jacobi)=\tel{3}\.\tr(\jacobi\.\const^T\widetilde{\sigma}^T)\,,
\end{align*}
thus in the cases $(\alpha)$ and $(\beta)$,
\begin{align*}
	\hpressure&=\tel{3}\.\tr(G\inv\widetilde{\sigma}^T)=\tel{3}\.\tr(\widetilde{\sigma}\.G\inv)\tag{5.5$\alpha$}\\
	\text{and}\qquad\qquad&\\
	\hpressure&=\tel{3}\.\tr \widetilde{\sigma}^T=\tel{3}\.\tr \widetilde{\sigma}\,.\tag{5.5$\beta$}
\end{align*}
Again, the mixed-variant definition $(\beta)$ yields the simpler formula.
\par\medskip\noindent
The infinitesimal deformation $\id+\intd{F}$ corresponds to the deformation $\id+\intd{\deform}$ in arbitrary coordinates according to: $\jacobi\.(\id+\intd{F})\.\intd{\force_0}=(\id+\intd{\deform})\.\jacobi\.\intd{\force_0}$. Thus $\intd{F}=\jacobi\inv\intd{\deform}\.\jacobi$
and, due to \eqref{5.2} and \eqref{5.4},
\begin{align*}
	\intd{\work}=\tr(\const^T\widetilde{\sigma}^T\jacobi\.\jacobi\inv\intd{\deform}\.\jacobi)=\tr(\jacobi\.\const\.\widetilde{\sigma}\.\intd{\deform})\,.	
\end{align*}
In the cases $(\alpha)$ and $(\beta)$ we find
\begin{align*}
	\intd{\work}&=\tr(G\inv\widetilde{\sigma}^T\intd{\deform})=\tr(\intd{\deform}^T\widetilde{\sigma}\.G\inv)\tag{5.4$\alpha$}\\
\text{and}\qquad\qquad&\\
	\intd{\work}&=\tr(\widetilde{\sigma}^T\intd{\deform})=\tr(\intd{\deform}^T\widetilde{\sigma})\,.\tag{5.4$\beta$}
\end{align*}
Again, we obtain a simpler result for the definition $(\beta)$.	

\subsection{Invariance of the law of elasticity}

For isotropic materials in Cartesian coordinates the law of elasticity has the form\footnote{See H.\ Richter: \enquote{Das isotrope Elastizitätsgesetz}, \emph{Zeitschrift für Angewandte Mathematik und Mechanik} 28.7/8 (1948), Pp.\ 205-–209 \cite{richter1948}.}
\begin{align*}
	e^j\.\sigma=\pdd{E}{j}\.\id+2\.\pdd{E}{k}\.L+3\.\pdd{E}{l}\.L^2\qquad\text{for}\qquad j=\tr L\,,\enskip k=\tr L^2\enskip\text{and}\enskip l=\tr L^3\,,
\end{align*}
where $E$ is the elastic potential per unit volume of the initial state.\transcomment{$e^j\sigma=\det F\cdot\sigma=\tau$ is the Kirchhoff stress tensor.}
\par\medskip\noindent
If we want this simple form to hold for arbitrary coordinates as well, then $\widetilde{\sigma}$ and $L$ must have the same mixed invariance, since the invariants and functional dependences are transferred only in this case. Therefore, and due to reasons mentioned above, it appears most practical to define both $\widetilde{\sigma}$ and $E$ covariant-contravariant, which is the reason this variance has been emphasized in the definition of $E$ in chapter 3.\\[1.4cm]
Received 25.~May, 1948\newpage
\section*{Review by Ruth Moufang (Zentralblatt für Mathematik und ihre Grenzgebiete)}
Hencky introduced the logarithms of the principal strains as quantities of strain for finite deformation of isotropic materials. Here, this definition of the strain tensor is recovered as a special case of a characterization based on the following postulates, where $F$ is the matrix of the linear transformation of the coordinate differentials and $G$ is the fundamental tensor of the metric:
\begin{enumerate}
	\item
		The strain tensor $E(F)$ is determined by the matrix $F$ and, apart from $F$, only depends on $G$.
	\item
		If $R$ is a rotation, then $E(F\.R)=E(F)$.
	\item
		A superposition principle holds such that for two coaxial stretches $V_1$ and $V_2$ and the corresponding strain tensors $E_1=E(V_1)$, $E_2=E(V_2)$, $E_3=E(V_1\.V_2)$, there exists a uniquely invertible function $f(x)$ with $f(E_1)+f(E_2)=f(E_3)$.
	\item
		For infinitesimal deformations $\id+\intd{\deform}$ in Cartesian coordinates, the strain tensor turns into\linebreak $\frac12\.(\intd{\deform}+\intd{\deform}^T)+o(\intd{\deform})$, where $o(x)$ denotes the usual symbol and $\deformT$ denotes the transpose of the matrix $\deform$ in general.
\end{enumerate}
If $F$, in Cartesian coordinates, is split into a product of a pure stretch $V$ with $3$ real positive eigenvalues and a Euclidean transformation, then the above postulates yield $E=f\inv(\log V)$, where $f\inv(x)$ is the inverse function of $f(x)$ and attains the form $x+o(x)$ for small $x$. In the simplest case one has to set $f\equiv x\equiv f\inv$, which leads to Hencky's approach. Moving to curvilinear coordinates then yields a covariant, contravariant or mixed tensor at choice. In the latter case, $E=f\inv(L^\ast)$ with $L^\ast=\frac12\log(G\.F\.\widehat{G}\inv F^T)$, where $\widehat{G}$ is the fundamental tensor with respect to the end position. Here, in general, both $f$ and $f\inv$ are tensor-valued functions of a tensor, e.g.\ given in the form of a convergent infinite series with a tensorial argument. --- Then the logarithm of the volume dilation is given by $\tr f(E)$, i.e.\ the trace of $f(E)$. --- The otherwise common strain tensor introduced by Trefftz\transcomment{\enquote{Trefftz's strain tensor} is the \enquote{Almansi strain tensor} $\frac{1}{2}(\id-B\inv)$ in the current configuration.} satisfies the above postulates for the superposition function $f(x)= -\frac12\log(1-2\.x)$. --- The strain deviator $D$ is deduced from the strain tensor by the requirements that two deformations which differ only by a similarity transformation have the same deviator and that the tensor of a volume preserving deformation is equal to its deviator. If, in general, the common deviator operation with respect to $E$ is denoted by $\dev E$, then $D=f\inv(\dev L^\ast)$. The discussion of the characteristic equation corresponding to $\dev L$ gives some indication of the physical meaning of the relation between $\tr(\dev L^3)^2$ and $\tr(\dev L^2)^3$ and indicates that $\sqrt{\tr(\dev L^2})$ can generally be considered a measure for the change of shape in agreement with the usual definition for infinitesimal deformations. --- The author refers the stresses to the undeformed surface element and defines the stress tensor via the requirements that
\begin{enumerate}
	\item
		in Cartesian coordinates, the force $\intd{A_0}$ acting on a surface element $\intd{\force_0}$ is given by $\intd{\force_0}=\sigma\.\intd{A_0}$,
	\item
		in curvilinear coordinates, $\widetilde{\sigma}$ is a tensor (or a tensor density),
	\item
		translating the surface element by $\intd{z}$ corresponds to the work $\intd{\work}=\iprod{\intd{z},\widetilde{\sigma}\.\intd{A}}$. 
\end{enumerate}
These conditions yield a representation of $\widetilde{\sigma}$ in terms of $\sigma$ as a mixed or twice-contravariant tensor. However, in the former case, $\widetilde{\sigma}$ is no longer symmetric along with $\sigma$. --- Computing the power for infinitesimal strain yields the known formulae and shows the advantage of using mixed tensors.
\par\medskip\noindent
\begin{flushright} 
	Ruth Moufang (Frankfurt a. M., 1950)
\end{flushright}
\section*{Review by William Prager (Mathscinet)}

To define strain in a continuous medium which undergoes a finite deformation, the author starts with the matrix $F$ which represents the mapping of a neighborhood of a point $\widehat{x}$ in the undeformed medium on to a neighborhood of the corresponding point $x$ in the deformed medium: $dx=F\,d\widehat{x}$. In a plastic material the history of deformation is important and, hence, the knowledge of $F$ alone is not sufficient. For an elastic material, on the other hand, $F$ completely characterizes the deformation. For an anisotropic elastic material, the rigid body rotation contained in $F$ is important, and $F$ itself must be used to describe the deformation. For an isotropic elastic material, however, this rigid body rotation is unessential; the strain tensor is then obtained by eliminating this rigid body rotation in a suitable manner. The author proceeds to establish postulates which should be satisfied by any acceptable definition of the strain tensor $E$. First of all, it must be possible to build up this tensor from the elements of the matrix $F$. Secondly, the tensor should not be influenced by a rigid body rotation which precedes the deformation characterized by the matrix $F$. Thirdly, if $V_1$ and $V_2$ denote pure stretches with coincident principal axes, $E_1=E(V_1)$ and $E_2=E(V_2)$ the corresponding strain tensors and $E=E(V_1\,V_2)$ the strain tensor corresponding to the deformation characterized by $V_1\,V_2$(=$V_2\,V_1$), there should exist a monotonic function $f(E)$ such that $f(E_1)+f(E_2)=f(E)$. Finally, the definition of the strain tensor should reduce to the customary one when infinitesimal deformations are considered. The author introduces a logarithmic strain tensor and shows that it satisfies these postulates.
\par\medskip\noindent
\begin{flushright} 
	William Prager (1949)
\end{flushright}
\end{refsegment}
\section*{Footnote by C. Truesdell and R. Toupin}

Later [1949] Richter worked out various special properties of [$\log V$] and [$\log U$]. Noticing that the condition of vanishing in uniform dilation does not determine a unique strain measure, Richter proposed a set of \underline{axioms}, including a \underline{superposition principle} for coaxial stretches, and showed that there are at $x$ and $X$ unique distortion tensors which satisfy them. This corrects an early attempt by Moufang \cite{moufang1947volumtreue}. Richter's distortion tensors are complicated algebraic functions of $e$ and $E$, respectively.
\par\medskip\noindent
\begin{flushright} 
	Clifford Truesdell and Richard Toupin (1960) \cite[p.270]{truesdell60}
\end{flushright}
\section*{Footnote by C. Truesdell and W. Noll}
The first attempts at mathematical treatment of Cauchy's idea [of an elastic material], apparently, are those of Reiner \cite{reiner1948elasticity}, Richter \cite{richter1948} and Gleyzal \cite{gleyzal1949};\\
Richter \cite{richter1952elastizitatstheorie} was the first to observe that the reduction follows at once from a simple and natural requirement of invariance, which is in fact a special case of the principle of material frame-indifference.
\begin{flushright} 
	Clifford Truesdell and Walter Noll (1965) \cite[p.119]{truesdell65}
\end{flushright}%

\section*{List of Symbols}
\begin{table}[H]
\vspace*{-2.45em}
\small
\begin{longtable}{| l | l | l |}
\hline
\vphantom{$\displaystyle\int$}{\normalsize\textbf{Our notation}}	&{\normalsize\textbf{Richter's notation}}				&\\
\hline
&&\\[-.49em]
$A$, $B$								&$\textfrak{A}$, $\textfrak{B}$/$\textfrak{C}$	&arbitrary $3\times3$-matrices\\
$a_{ik}$, $(A)_{ik}$	&$a_{ik}$, $(\textfrak{A})_{ik}$						&entry in the $i$-th row and the $k$-th column of $A$\\
$\det A$								&$|\textfrak{A}|$								&determinant of $A$\\
$\tr A$									&$\{\textfrak{A}\}$								&trace of $A$\\
$A^T$									&$\overline{\textfrak{A}}$						&transpose of $A$\\
$\id$									&$\textfrak{E}$									&identity tensor\\
$A\inv$								&$\textfrak{A}\inv$							&inverse of $A$\\
$x$, $y$, \dots							&$\textfrak{x}$, $\textfrak{y}$, \dots			&vectors\\
&&\\
$F$, $\deform$						&$\textfrak{A}$, $\textfrak{B}$		&Jacobian matrices (deformation gradients)\\
$\widehat{x}$							&$\widehat{x}$							&preimage of $x$ under $F$\\
$E(F)$, $E$						&$\textfrak{B}(\textfrak{A})$, $\textfrak{B}$	&strain tensor corresponding to $F$\\
$R$										&$\textfrak{R}$							&pure Euclidean rotation\\
$V$										&$\textfrak{S}$							&pure stretch\\
$\widehat{\hphantom{A}}$				&$\widehat{\hphantom{A}}$				&indicator of a tensor being associated with the reference configuration $\widehat{x}$\\
$V_1$, $V_2$							&$\textfrak{S}_1$, $\textfrak{S}_2$		&coaxial stretches\\
$E_1$, $E_2$							&$\textfrak{V}_1$, $\textfrak{V}_2$		&strain tensors $E(V_1)$, $E(V_2)$\\
$f$										&$f$									&uniquely invertible function with $f(E_1)+f(E_2)=f(E)$\\
$o$										&$o$									&function with $y=o(x)$:\enspace $\lim \frac{y}{x}=0$\\
&&\\
$\opoint$, $\ipoint$ 				&$\textfrak{h}$, $\textfrak{y}$ 		&original point and its image under the deformation $\deform$\\
$E_0$ 									&$\textfrak{W}$ 						&strain tensor with respect to $\deform$\\
$Z$										&$\textfrak{Z}$							&$Z=f(E_0)$\\
$L$										&$\textfrak{L}$							&logarithmic strain tensor:\enspace $L=\log V$\\
$f\inv$								&$g$									&inverse function of $f$\\
$h$, $k$								&$h$, $k$								&functions:\enspace $h(x)=f\inv(\log(x))$, $k(x)=h(\sqrt{x})$\\
&&\\
$\jacobi$								&$\textfrak{U}$							&Jacobian matrix of $x=x(q)$\\
$G$										&$\textfrak{G}$							&metric fundamental tensor\\
$L^{\ast}$								&$\textfrak{L}^{\ast}$					&logarithmic strain tensor in curvilinear coordinates\\
$v$										&$v$									&dilatation being associated with $F$:\enspace $v=\det F$\\										
$T$										&$\textfrak{T}$							&\enquote{common} strain tensor, $T=\frac{1}{2}(\id-B\inv)$, Almansi strain tensor \\
&&\\
$D$										&$\textfrak{D}$							&strain deviator (change of shape)\\												
$\dev A$								&$\widetilde{\textfrak{Q}}$				&common deviator of the matrix $A$: $\dev A=A-\frac{1}{3}\tr (A)\cdot\id$\\
$\zeta$									&$\zeta$								&$\zeta$ characterizes the kind of loading\\
&&\\
$\sigma$								&$\textfrak{P}_0$						&Cauchy stress tensor\\
$\widetilde{\sigma}$					&$\textfrak{P}$							&stress tensor in curvilinear coordinates\\
$\intd{A}$								&$\intd{\textfrak{f}}$					&surface element\\
$\intd{\force}_0$						&$\intd{\textfrak{k}}_0$				&the force acting on $\intd{A}_0$ at $\ipoint$\\				
$\const$								&$\textfrak{C}$							&constant\\
$\intd{\work}$							&$\intd{A}$								&differential of the expended work\\
$V$										&$V$									&volume\\
$\mathcal{F}$							&$F$									&surface of $V$\\
$\hpressure=\frac{1}{3}\tr \sigma$		&$\sigma$								&hydrostatic stress, mean stress\\
$\iprod{x,y}$							&$x\cdot y$								&scalar product\\
$|x|^2$									&$x^2$									&squared length of a vector\\[-.49em]
&&\\
\hline			
\end{longtable}
\caption{\label{table:notation}Changes made to Richter's notation.}
\end{table}
\renewcommand{\url}[1]{\href{#1}{www.uni-due.de}}
\printbibliography
\end{document}